\documentclass{ifacconf}

\usepackage{graphicx}      
\usepackage{natbib}        
\usepackage{amsmath}
\usepackage{amssymb}
\usepackage{color}
\usepackage{epstopdf}
\usepackage{algorithm}
\usepackage[noend]{algpseudocode}
\usepackage{mathtools}

\newcommand{\Vect}{\mathrm{vec}}

\makeatletter
\def\BState{\State\hskip-\ALG@thistlm}
\makeatother
\newcommand*\xpartial[4]{%
  \mathrlap{\frac{\partial #1_{\mathrlap{#3}}}{\partial #2_{\mathrlap{#4}}}}%
  \phantom{\frac{\partial #1_{#3}}{\partial #2_{#4}}}%
  }
\begin{document}
\begin{frontmatter}
\hyphenation{Crass-idis}
\title{Large-Scale Distributed Kalman Filtering via an Optimization Approach} 

\thanks[footnoteinfo]{This research was supported by the U.S. Army Research
Laboratory and the U.S. Army Research Office under contract number W911NF-13-1-0340 and AFOSR grant  FA9550-16-1-0022.}

\author{Mathias Hudoba de Badyn$^*$ Mehran Mesbahi$^*$}

\address{$^*$William E. Boeing Department of Aeronautics and Astronautics, University of Washington, Seattle, WA 98195 USA \\(e-mails: \{hudomath,mesbahi\}@uw.edu).}

\begin{abstract}                
Large-scale distributed systems such as sensor networks, often need to achieve filtering and consensus on an estimated parameter from high-dimensional measurements. 
Running a Kalman filter on every node in such a network is computationally intensive; in particular the matrix inversion in the Kalman gain update step is expensive.
In this paper, we extend previous results in distributed Kalman filtering and large-scale machine learning to propose a gradient descent step for updating an estimate of the error covariance matrix; this is then embedded and analyzed in the context of distributed Kalman filtering.
We provide properties of the resulting filters, in addition to a number of applications throughout the paper.
\end{abstract}

\begin{keyword}
Machine learning, fast Kalman algorithms, state estimation, gradient methods
\end{keyword}

\end{frontmatter}

\section{Introduction}

The Kalman filter is an algorithm that uses the known dynamics of a system to remove noise from measurements of that system.
When considering large-scale dynamical systems, implementation of the standard or extended Kalman filters can be computationally difficult.
In such cases, the Kalman filter requires the inversion of very large matrices at each timestep.
This may cause the Kalman filter to run slower than the dynamical process it is trying to measure, or to severely reduce the temporal resolution of the measurements.

There are many systems for which measurements are taken by a network of sensors and are of high dimension; examples of such systems can be found in \cite{khan2008} and \cite{Kutz2013}.
Previous methods for circumventing this problem include decomposing the dynamical system being measured into several subsystems and distributing the subsystems over the sensor network as in~\cite{khan2008}, or using Monte-Carlo methods for estimating the error covariance, such as in~\cite{Furrer2007}.

Previously,~\cite{Sutton1992} proposed to modify the Kalman filter error covariance update with a gradient descent method for the purpose of mimimizing memory consumption, albeit for a specific instantiation of a SISO linear system. 
In this paper, we extend the gradient descent algorithm for estimating the Kalman filter error covariance to the general MIMO linear system as a proposed solution to the problem of running a Kalman filter on a high-dimensional system.
We improve the gradient descent using Nesterov acceleration and adaptive learning rate methods.
Lastly, we apply the methods above to distributed Kalman filtering on a sensor network.

Distributed Kalman filtering seeks to estimate the state of a system by distributing the tasks of measuring the system and subsequently filtering the data to many agents, who then collectively assemble the state estimate.
Such algorithms utilizing consensus to provide a global estimate were presented by~\cite{Olfati-Saber2005} and~\cite{Olfati-Saber2007}, and later extended by~\cite{Carli2008}.
Performance of distributed Kalman filters using graph-theoretic quantities were studied by~\cite{Spanos2005}.

The organization of the paper is as follows. In \S\ref{sec.math}, we outline the mathematical notation and conventions we use.
We discuss the gradient descent algorithm, Nesterov acceleration and adaptive learning rate methods in \S\ref{sec.filt}, and the distributed version of these algorithms with a relevant example in \S\ref{sec.dist}.
The paper is summarized in \S\ref{sec.conc}.
\section{MATHEMATICAL PRELIMINARIES}\label{sec.math}
In this section, we lay out the mathematical notation used in this paper, and summarize the essential background on Kalman filtering.

Consider the \emph{noisy discrete-time controlled linear system} 
\begin{align*}
  \begin{aligned}
     x_{k+1} &= A x_k + B u_k + \Upsilon _k w_k\\
\tilde  {y}_k &= C_k x_k + v_k,
  \end{aligned}\label{eq.linsys}
\end{align*}
where $x_k\in\mathbb{R}^n$ denotes the \emph{state vector}, $u_k\in\mathbb{R}^m$ denotes the \emph{control vector}, $\tilde y_k\in\mathbb{R}^p$ is the \emph{output vector} or \emph{measurement}, and $w_k,v_k$ are Gaussian white noise vectors of appropriate dimensions with covariances $Q_k$ and $R_k$ respectively.
The subscript $k$ refers to the timestep.
The random vectors $w_k$ and $v_k$ represent system disturbances and sensor noise respectively, and the matrix $\Upsilon_k$ describes how the disturbance propagates into the system.
Since most matrix quantities will have a subscript $k$ denoting the timestep, we will refer to the $ij$th entry of a matrix $A_k$ with bracketed superscripts: $A_k^{(ij)}$.
Similarly, the $i$th entry of a vector $a_k$ is denoted $a_k^{(i)}$, the $i$th column of a matrix $A_k$ is denoted $A_k^{[i]}$ and the $i$th row is denoted $A_k^{\{i\}}$.

A network is represented by a \emph{graph} $\mathcal{G}=(V,E)$ where $V$ is a set of \emph{nodes} representing agents in the network, and $E$ is the set of edges $\{i,j\}$ representing the connections between agents $i$ and $j$.
The \emph{neighbourhood set} $N_i$ of node $i$ is the set of indices $j$ such that $\{i,j\}$ is an edge.
We assume for simplicity that the graph is undirected, meaning that $\{i,j\}$ and $\{j,i\}$ represent the same edge.

The purpose of Kalman filtering is to provide an accurate estimate $\hat{x}_k$ of the state $x_k$ using the measurement $\tilde y_k$ and the known information about the system and noise.
The standard \emph{discrete-time linear Kalman filter} yields the estimate $\hat{x}_{k+1}$ and is given by 
\begin{align*}
 \hat x_{k+1} &= A\left(\hat x_k + K_k \left[ \tilde y_k - C \hat{x}_k\right] \right)\\
  K_k &= P_k C^T \left[ C P_k C^T + R_k\right]^{-1}\\
 P_{k+1} &= A \left(I-K_kC\right)P_k A^T + \Upsilon_k Q_k \Upsilon_k^T,
\end{align*}
where $K_k$ is the \emph{Kalman gain} and $P_k$ is the \emph{error covariance matrix}.
The problem for high-dimensional systems is the inverse when computing $K_k$. 
In this paper, we assume that the sensor noise is uncorrelated, and so $R_k$ is a diagonal positive-definite matrix.
This leaves the term $C P_k C^T$; if this term is diagonal, then the matrix inverse becomes a series of $n$ scalar divisions.
We propose to replace $P_k$ with a diagonal estimate $\hat{P}_k$, and one can also assume that $C$ is a matrix that measures individual states of the system without redundancy.
Therefore, a system with $1\leq p \leq n$ outputs will have a measurement matrix of the form $C = [\mathrm{diag}(C^{(ii)})I_{p}~\mathbf{0}_{n-p}]$.
This assumption is common for networked systems where one measures the state of individual nodes of the network, see for example \cite{Chapman2015a}.
In the next few sections, we will formalize these assertions in the context of our algorithm.
 
\section{GENERAL KALMAN FILTERING}\label{sec.filt}

In this section, we extend the gradient descent estimate for the error covariance proposed by~\cite{Sutton1992}.
In \S\ref{sec.gdd}, we discuss the linear Kalman filter, and then extend the error covariance estimate to an accelerated gradient descent in \S\ref{sec.acc}.
The accelerated gradient descent benefits from a clever adaptive learning rate, which is discussed in \S\ref{sec.alr}.
The final algorithm combining all these methods is summarized in Algorithm~\ref{alg.1}.

\subsection{Gradient Descent for the Error Covariance Update: Linear Case}\label{sec.gdd}

Gradient descent is an iterative algorithm that seeks to find the local minimum of a function $f$ by stepping in the direction of the largest gradient: 
\begin{align*}
  \beta_{k+1} = \beta_k - \mu \nabla f(\beta_k),
\end{align*}
where $\mu$ is the \emph{learning rate}.
For the Kalman filter, we seek to find a diagonal matrix estimate $\hat{P}_k$ of $P_k$.
We do this by assuming $\hat{P}_k$ is of the form $\hat{P}_k = \mathrm{diag}(\hat{P}_k^{(ii)})=\mathrm{diag}(e^{\beta_k^{i}})$,
where $\beta_k$ is a parameter undergoing gradient descent attempting to minimize the norm of the error $\delta_k = \tilde{y}_k-C\hat{x}_k$: 
\begin{align*}
    \beta_{k+1}^{(i)} = \beta_k^{(i)} - \frac{1}{2}\mu \frac{\partial(\delta_k^T\delta_k)}{\partial \beta_k^{(i)}}.
\end{align*}
Computing the gradient of $\delta_k^T\delta_k$ as outlined in Appendix~\ref{app.lin} yields the following gradient descent equations for $\hat{P}$:
\begin{align*}
  \beta_{k+1}^{(i)} &= \beta_k^{(i)} + \mu \delta_k^T C^{[i]} h_k^{(i)}\\
  h_{k+1}^{(i)} &= h_k^{(i)}\left(1-k_k^{(i)T}C^{[i]}\right)^+ + \left(k_k^{(i)} - k_k^{(i)}k_k^{(i)T}C^{[i]}\right)^T \delta_k
\end{align*}
\begin{equation}
\hat{P}_k^+ = \mathrm{diag}\left(\exp\left(\beta_{k+1}^{(i)}\right)\right)\label{eq:1}~~~~~~~~~~~~~~~~~~~~~~~~~~~~~~~~ 
\end{equation}
The form of $\hat{P}^+$ guarantees that it remains positive-definite, which is required to preserve convergence of the Kalman filter.

\subsection{Nesterov-Accelerated Methods}\label{sec.acc}
Nexterov acceleration is a method used to increase the convergence rate of gradient descent [\cite{nesterov1983}].
Although Nesterov-accelerated gradient descent converges in fewer timesteps, it does so by sacrificing monotonicity: the gradient descent trajectory will tend to oscillate as it converges towards the estimate.

In order to implement Nesterov acceleration, one must see that the quantity $h_k^{(i)}$ is a function of $K_{k-1}$, which in turn is a function of $\beta^{(i)}_k$: 
\begin{align*}
  \beta_{k+1}^{(i)} &= \beta_k^{(i)} - \dfrac{1}{2}\mu \nabla f\left(\beta_k^{(i)}\right)\\
&= \beta_k^{(i)} + \mu \delta_k ^T C^{[i]}h_k^{(i)} \left( K_{k-1}\left(\beta_k^{(i)}\right)\right).
\end{align*}
Hence, we can get the form of the Nesterov-accelerated gradient descent by introducing the quantity $\alpha_k$:
\begin{align*}
  \beta_{k+1}^{(i)} &= \alpha_k^{(i)} - \dfrac{1}{2} \mu \nabla f \left( \alpha_k^{(i)}\right)\\
&= \alpha_k^{(i)} +\mu\delta_k^T C^{[i]} h_k^{(i)} \left( K_{k-1}\left( \alpha_k^{(i)}\right)\right)
\end{align*}
where $\alpha_k^{(i)}$ is given by 
\begin{align*}
  \alpha_k^{(i)} = \beta_k^{(i)} + \dfrac{k-1}{k+2} \left(\beta_k^{(i)} - \beta_{k-1}^{(i)}\right).
\end{align*}
Therefore, the Nesterov-accelerated estimate $\hat P$ replaces the covariance update with the following equations:
\begin{align*}
\hat P &=  \mathrm{diag}\left( e ^{\beta^{(i)}_{k+1}}\right),~k_k^{(i)} = \exp\left({\alpha_{k+1}^{(i)}}\right) D^{-1}_k C^{[i]}
\end{align*}
\begin{align}
  \beta_{k+1}^{(i)} &= \alpha_k^{(i)} +\mu\delta_k^T C^{[i]} h_k^{(i)} \left( K_{k-1}\left( \alpha_k^{(i)}\right)\right), \label{eq.nesterov}
\end{align}
where $h_k$ is as in Equation~\eqref{eq:1}, but using this new expression for $k_k^{(i)}$ in Equation~\eqref{eq.nesterov}.
\subsection{Adaptive Learning Rate Methods}\label{sec.alr}

Using a constant learning rate $\mu$ can lead to suboptimal performance, and therefore it is prudent to adapt $\mu$ on-the-fly as the gradient descent is run.
A well-known method for adapting the learning rate is given by Barzilai and Borwein~\cite{Barzilai1988}.
Suppose we have the gradient descent with learning rate $\mu$:
\begin{align*}
  \beta_{k+1}^{(i)} = \beta_{k}^{(i)} - \frac{1}{2}\mu \nabla f\left(\beta_{k}^{(i)}\right).
\end{align*}
Let $\Delta\beta = \beta_k-\beta_{k+1}$ and let $\Delta g(\beta) = \nabla f \left(\beta_{k+1}\right)-\nabla f \left(\beta_{k}\right)$.
Then, the Barzilai and Borwein method selects the learning rate parameter using a secant line approximation:
\begin{align*}
  \frac{1}{2}\mu_k = \arg\min_{\lambda} \left\| \Delta \beta - \lambda \Delta g(\beta) \right\|.
\end{align*}
For accelerated gradient descent, the updated learning rate $\mu$ is given by
\begin{align*}
  \mu_k = 2\dfrac{\Delta g(\beta)^T\Delta \beta}{\Delta g(\beta)^T\Delta g(\beta)},
~  \Delta \beta = \beta_k - \alpha_{k-1}.
\end{align*}
Combining the accelerated gradient descent estimate of the error covariance with the adaptive learning rate $\mu_k$ yields the final algorithm as outlined in Algorithm 1.
\begin{algorithm}\label{alg.1}
\caption{Discrete Kalman Filter with Accelerated Gradient Descent on $P$ and Adaptive Learning Rate $\mu$}\label{alg.discac}
\begin{algorithmic}[1]
\BState \emph{Initialize}
\State $\hat x(t_1) = \hat x_0$
\State $P_0 = \mathbb{E}\left\{\tilde x(t_1) \tilde x(t_1)^T\right\}$
\BState \emph{Loop}:
\For{$k=1$ to $k=t_f$}
\BState \emph{Gain}
\State $K_k = \hat{P}_k^-C^T(C\hat{P}_k^-C^T+R_k)^{-1}$
\BState \emph{Accelerated Gradient Descent and Adaptive $\mu_k$}
\State $\Delta \beta = \beta_k - \alpha_{k-1}$
\State $\Delta g(\beta) = \nabla f \left(\beta_{k+1}\right)-\nabla f \left(\beta_{k}\right)$
\State $D_k = R_k + C \hat P_k C^T$
\State  $\mu_k = 2(\Delta g(\beta)^T\Delta \beta)/(\Delta g(\beta)^T\Delta g(\beta))$
\State $\alpha_k^{(i)} = \beta_k^{(i)} + \dfrac{k+1}{k+2}\left(\beta_k^{(i)} - \beta_{k-1}^{(i)}\right)$
\State $\kappa_k^{(i)} = \exp\left(\alpha_k^{(i)}\right)D_k^{-1}C^{[i]}$
\State $\beta_{k+1}^{(i)} = \beta_k^{(i)} + \mu_k \delta_k^T C^{[i]} h_k^{(i)}$
\BState $h_{k+1}^{(i)} = h_k^{(i)}(1-\kappa_k^{(i)T}C^{[i]}) + (\kappa_k^{(i)} - \kappa_k^{(i)}\kappa_k^{(i)T}C^{[i]})^T \delta_k$
\BState \emph{Update}
\State $\hat x_k^+ = \hat x_k^- + K_k \left[ \tilde y_k -C\hat x_k^-\right]$
\BState \emph{Propagate}
\State $\hat x_{k+1}^-=A_k \hat x_k^+ + B_k u_k$
\State $\hat{P}_{k+1}^- = \mathrm{diag}\left(\exp\left(\beta_{k+1}^{(i)}\right)\right)$
\EndFor
\end{algorithmic}
\end{algorithm} 
\subsection{Filter Properties}

In the standard Kalman filter, $P_k$ is computed from an algebraic Ricatti equation.
In order for this equation to have a positive-definite solution, a requirement for stable error dynamics, one must assume that $(A^T,C^T)$ is stabilizable and $(A,\sqrt{Q})$ is detectable (see~\cite{Dorato1994}).
Therefore, since our estimate of $P$ is positive definite by construction and not computed by an algebraic Ricatti equation, we do not need to have these assumptions.
The next theorem shows that all one needs for stability of the error dynamics is for the system matrix $A$ to be stable.
\begin{thm}
  Suppose $A$ is an asymptotically stable matrix, and so the eigenvalues of $A$ satisfy $|\lambda(A)|\leq 1$.
  Then, the error dynamics 
  \begin{align*}
    \tilde{x}_{k} = x_{k} - \hat{x}_{k}^-
  \end{align*}
are stable, where $\hat{x}_k^-$ is the estimate computed according to Algorithm 1.
\end{thm}
\begin{pf}
  We can write the error dynamics as 
  \begin{align*}
    \tilde{x}_{k+1} &= x_{k+1} - \hat{x}_{k+1}^-\\
    &= (A-LC)\hat{x}_k^- + ( w_k - L v_k),
  \end{align*}
where $L$ is given by 
\begin{align*}
  L = A\hat{P}C^T(C\hat{P}C^T+R)^{-1}.
\end{align*}
Hence, it suffices to show that $(A-LC)$ has eigenvalues in the unit disc.
We can write $(A-LC)$ as 
\begin{align*}
 A-LC &= A- A\hat{P}C^T(C\hat{P}C^T+R)^{-1}C = A(I-N),
\end{align*}
where $N=\hat{P}C^T(C\hat{P}C^T+R)^{-1}C$.
By the form of $C$ and $\hat{P}$, if the system has $p$ outputs, then $N$ is a diagonal matrix with the form 
\begin{align*}
  N^{(ii)} =
  \begin{cases}
    \frac{e^{\beta_k^{(i)}}(C^{(ii)})^2}{e^{\beta_k^{(i)}}(C^{(ii)})^2 + R^{(ii)}}&\text{if }1\leq i \leq p\\
    0 & \text{if }p<i\leq n.
  \end{cases}
\end{align*}
It is clear that $0\leq N^{(ii)} <1$ since $R$ is a positive definite diagonal matrix.
Therefore, $0<(I-N)^{(ii)}\leq1$, and so it follows that 
\begin{align*}
 |\lambda(L)| &= |\lambda(A(I-N))|\\ &
\leq |\lambda_{\max}(A)||\lambda_{\max}(I-N)| \\ 
&\leq |\lambda_{\max}(A)|,
\end{align*}
and so by the asymptotic stability of $A$, the asymptotic stability of $(A-LC)$ follows.
\end{pf}
\subsection{Example}

\begin{figure}
  \centering
  \includegraphics[width=\columnwidth]{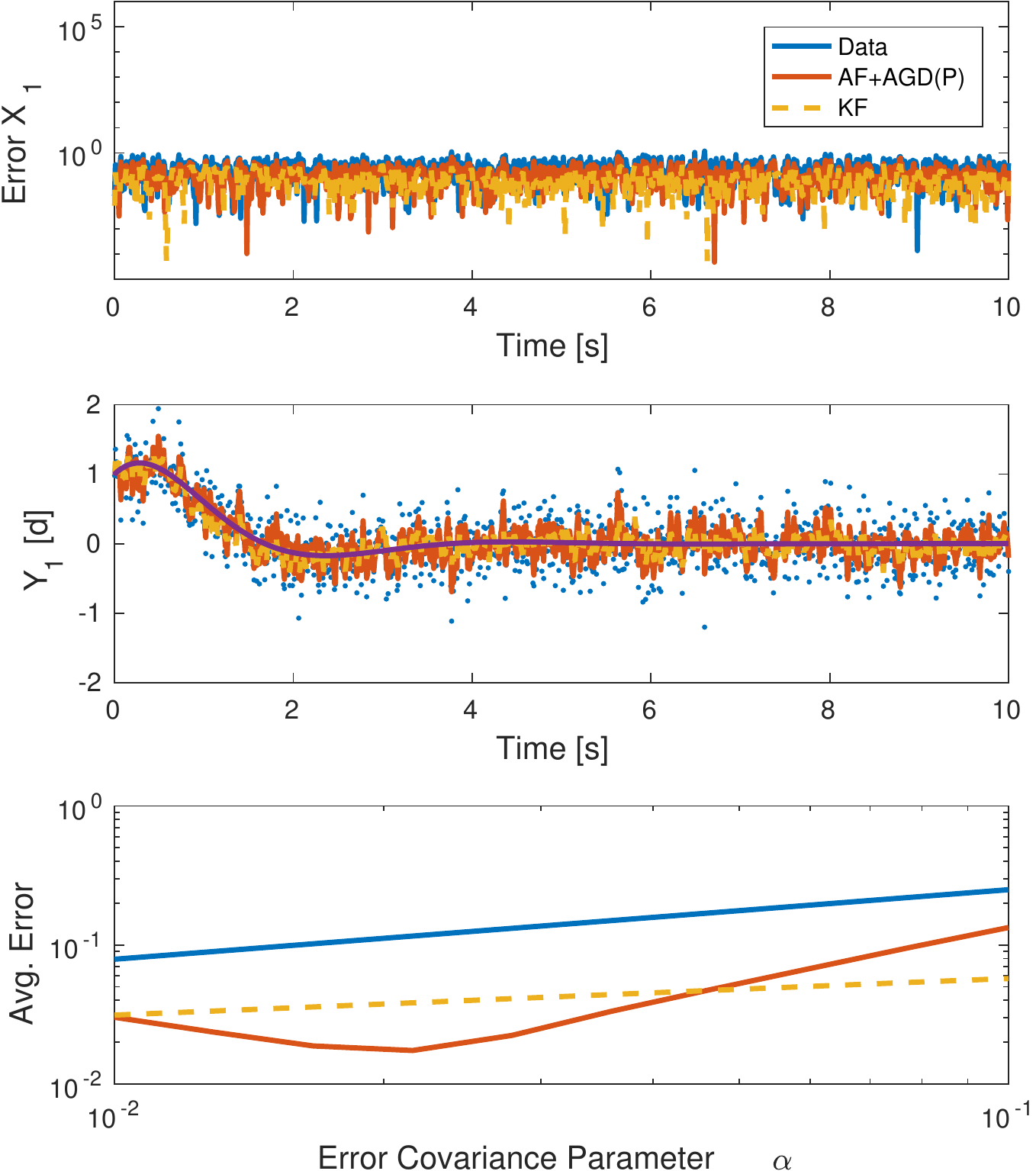}
  \caption{Numerical test of Algorithm~1. Top: Error from true state from data, accelerated filter and Kalman filter. Middle: $y_1=x_1$ for data, accelerated filter and Kalman filter. Bottom: Average steady-state error. }
  \label{fig.numeg}
\end{figure}

Here, we show a simple numerical implementation of Algorithm 1 to see how the algorithm performs compared to the standard Kalman filter for different strengths of noise.
Consider the linear system with matrices 
\begin{align*}
  A = \left[ 
  \begin{array}{cc}
    1.001&0.011 \\
    -0.0301 & 0.98
  \end{array}\right], B = \left[ 
              \begin{array}{c}
                5\mathrm{E}-5\\1\mathrm{E}-2
              \end{array}\right],C=I
\end{align*}
propagating with process noise of covariance $Q=0.001$ and measurement noise of covariance $R=\alpha I$ for various values of $\alpha$.
The first and second plots of Figure~\ref{fig.numeg} show the error and first state of the system for $R=0.115I$ for the raw measurement data, and the filtered data using the accelerated gradient descent and the standard Kalman filter.
The third plot shows the average steady state error computed on the time interval $t\in[6.6,10]$ for various values of $\alpha$ indicating the covariance of the measurement error.
Figure~\ref{fig.numeg} shows that the filter in Algorithm 1 achieves similar performance as the standard Kalman filter.

\section{DISTRIBUTED KALMAN FILTERING}\label{sec.dist}

\cite{Olfati-Saber2007} describes several methods of Kalman filtering using a network of sensors, with each $n$th sensor using a unique sensor model $\tilde y_{n,k} = C_{n,k} x_k + v_{n,k}$, to arrive at a global measurement of a system.
The context of their system model is on a graph $\mathcal{G}$, which is an abstract representation of the connection structure of sensors, which are represented by nodes.
An edge in the graph denotes a communication link between sensors by which they can exchange information about their current state estimate.
One such information exchange algorithm is called \emph{consensus}, where each node in the graph changes its estimated value by continuously averaging with the estimates of its neighbours: 
\begin{align*}
  \dot{\hat{x}}_i = \sum_{j\in N_i}\left(\hat{x}_j - \hat{x}_i\right).
\end{align*}

One particular algorithm in~\cite{Olfati-Saber2007} utilizes a \emph{Kalman-consensus filter} in which every node of the network computes an estimate of the state of the system and then performs consensus with its neighbours on this estimate.
This allows sensors that only see a portion of the system to collect a global state estimate over time.

For high-dimensional systems, implementing this filter is problematic again because each sensor has to invert a large matrix at each timestep in order to compute the Kalman filter gain.
In Algorithm 2, we summarize the Kalman-consensus filter,  where every node performs an accelerated gradient descent with adaptive learning rate, and shares its estimate with its neighbours. 
In the following section, we provide an example of such a high-dimensional system.

\begin{algorithm}\label{alg.1}
\caption{Discrete Distributed Kalman-Consensus Filter $\mu$}\label{alg.discac}
\begin{algorithmic}[1]
\BState \emph{Initialize}
\State $\hat x(t_1) = \hat x_0$
\State $P_0 = \mathbb{E}\left\{\tilde x(t_1) \tilde x(t_1)^T\right\}$
\State $\epsilon = t_2-t_1$
\BState \emph{Loop}:
\For{$k=1$ to $k=t_f$}
\For{each node $j\in V$}
\BState \emph{Aggregate Data}
\State $S_j = {|N_i|}^{-1}\sum_{l\in N_j}C_i$
\State $y_j = {|N_i|}^{-1}\sum_{l\in N_j}C^T_i \tilde y_i$
\BState \emph{Accelerated Gradient Descent and Adaptive $\mu_k$}
\State $\Delta \beta = \beta_k - \alpha_{k-1}$
\State $\Delta g(\beta) = \nabla f \left(\beta_{k+1}\right)-\nabla f \left(\beta_{k}\right)$
\State $D_k = R_k + C \hat P_k C^T$
\State  $\mu_k = 2(\Delta g(\beta)^T\Delta \beta)/(\Delta g(\beta)^T\Delta g(\beta))$
\State $\alpha_k^{(i)} = \beta_k^{(i)} + \dfrac{k+1}{k+2}\left(\beta_k^{(i)} - \beta_{k-1}^{(i)}\right)$
\State $\kappa_k^{(i)} = \exp\left(\alpha_k^{(i)}\right)D_k^{-1}C^{[i]}$
\State $\beta_{k+1}^{(i)} = \beta_k^{(i)} + \mu_k \delta_k^T C_n^{[i]} h_k^{(i)}$
\BState $h_{k+1}^{(i)} = h_k^{(i)}(1-\kappa_k^{(i)T}C_n^{[i]}) + (\kappa_k^{(i)} - \kappa_k^{(i)}\kappa_k^{(i)T}C_n^{[i]})^T \delta_k$
\BState \emph{Update Kalman-Consensus Estimate}
\BState $\hat x_{j,k}^+ = \hat x_{j,k}^- + D^{-1}_k\left[ y_j -S_j\hat x_{j,k}^- + \epsilon\sum_{l\in N_j}(\hat x_{l,k}^--\hat x_{j,k}^-)\right]$
\BState \emph{Propagate}
\State $\hat x_{k+1}^-=A_k \hat x_k^+ + B_k u_k$
\State $\hat{P}_{k+1}^- = \mathrm{diag}\left(\exp\left(\beta_{k+1}^{(i)}\right)\right)$
\EndFor
\EndFor
\end{algorithmic}
\end{algorithm} 

\subsection{Example}
\begin{figure}
  \centering
  \includegraphics[width=0.75\columnwidth]{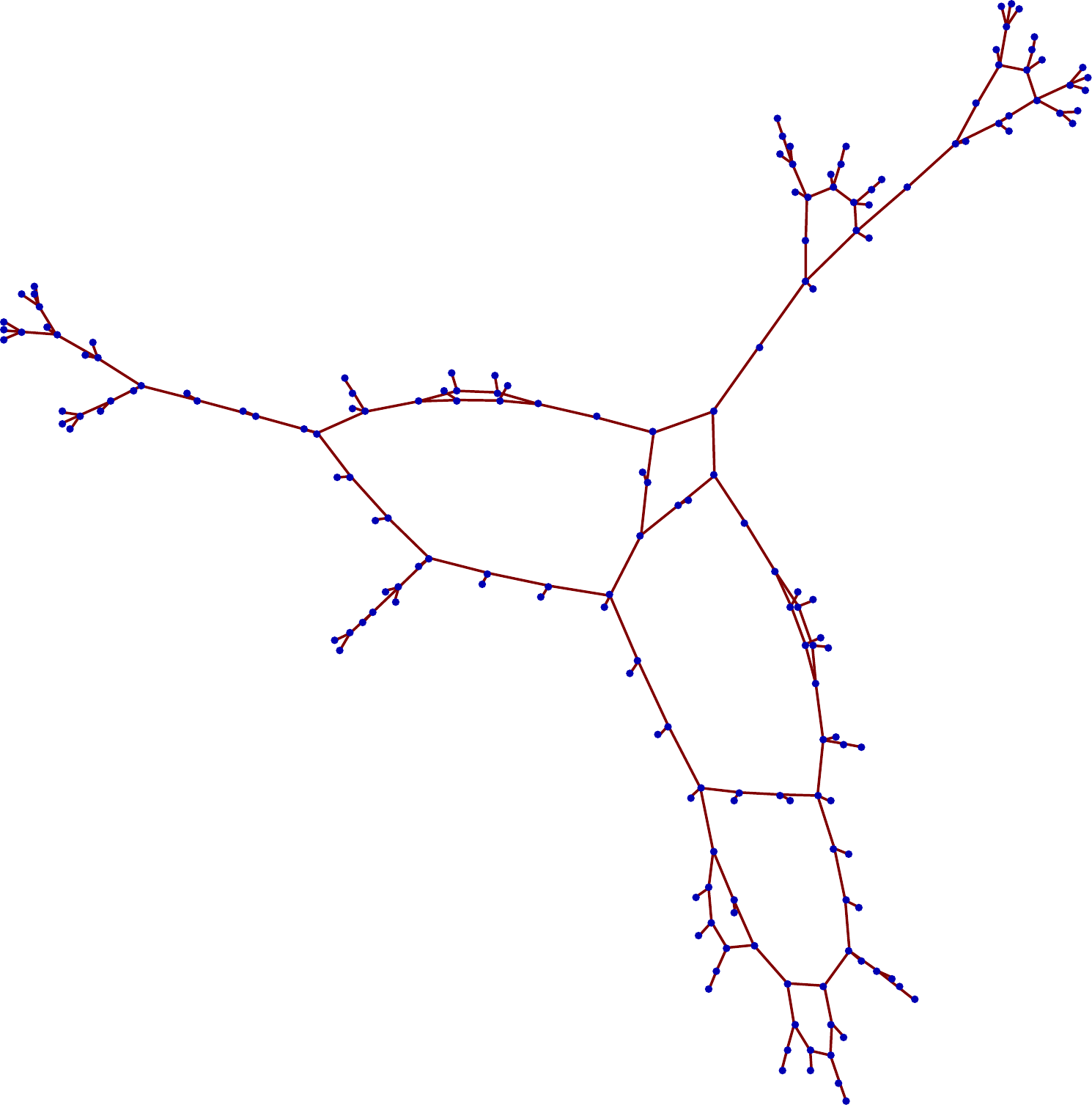}
  \caption{Sensor network graph with 178 nodes and 186 edges}
  \label{fig.graph}
\end{figure}
\begin{figure}
  \centering
  \includegraphics[width=\columnwidth]{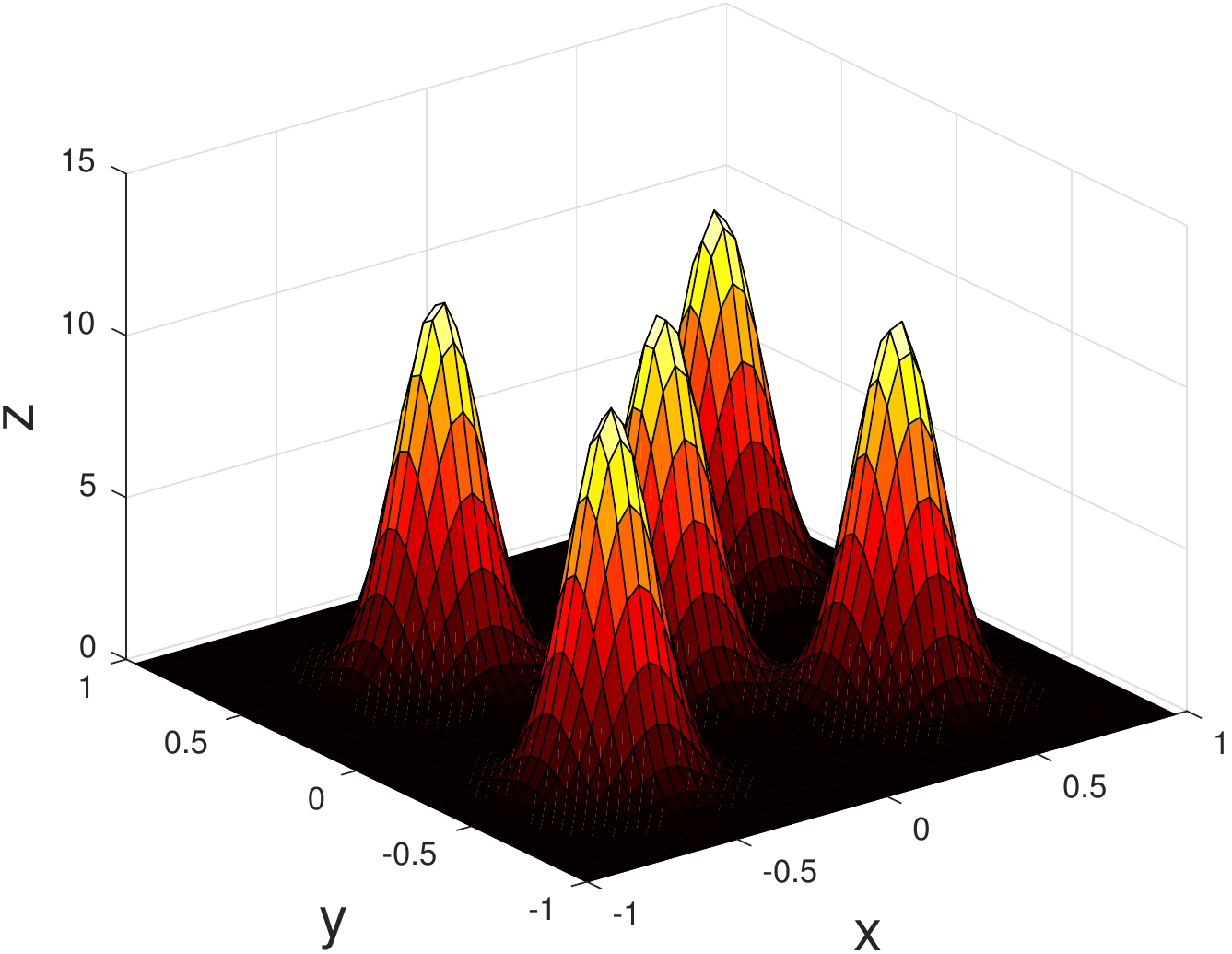}
  \caption{True state of system at $t=1.2$s}
  \label{fig.truestate}
\end{figure}
\begin{figure}
  \centering
  \includegraphics[width=\columnwidth]{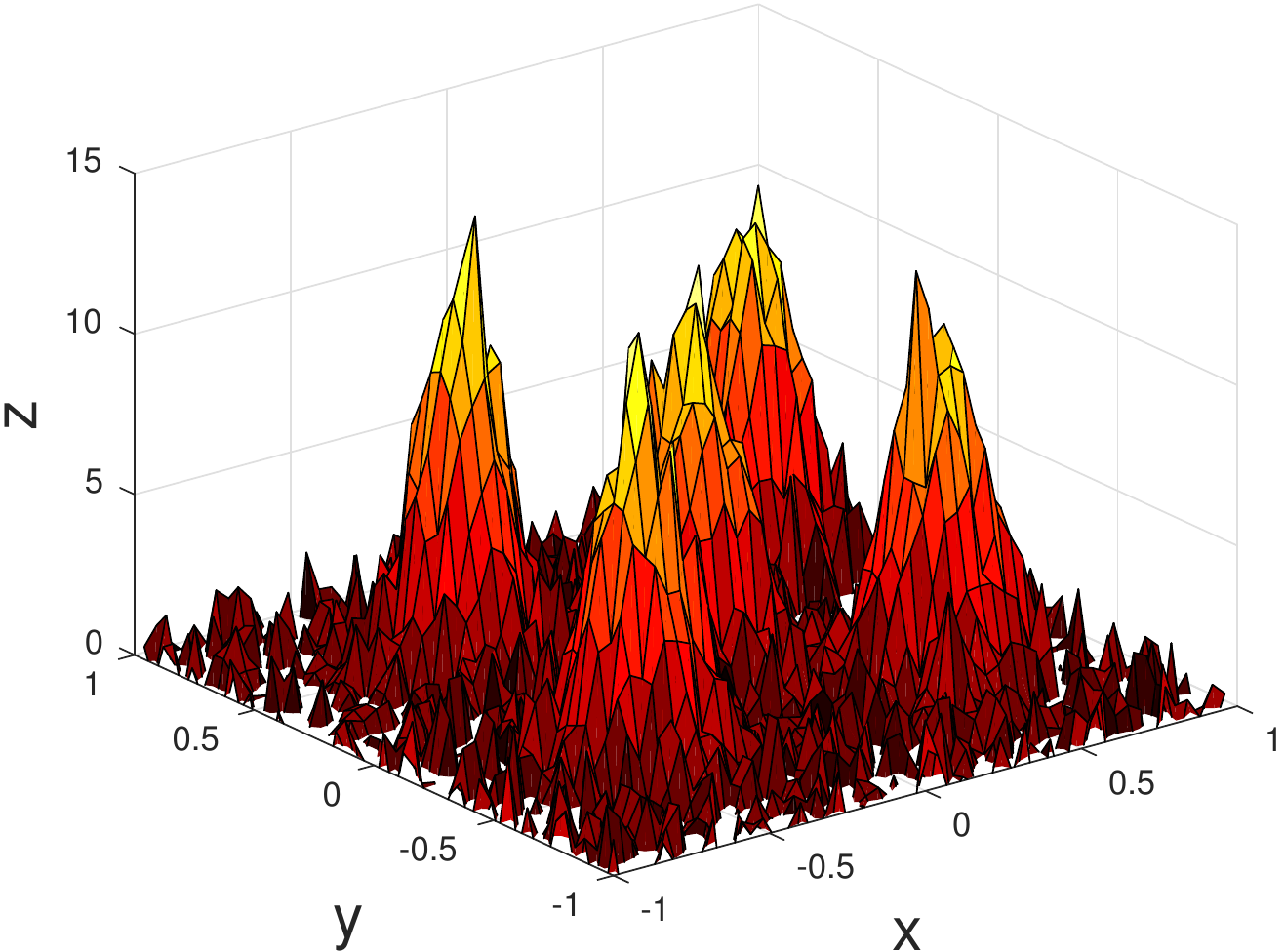}
  \caption{Sample measurement $\tilde y$ over the entire system at $t=1.2$s}
  \label{fig.ytilde}
\end{figure}
\begin{figure}
  \centering
   \includegraphics[width=\columnwidth]{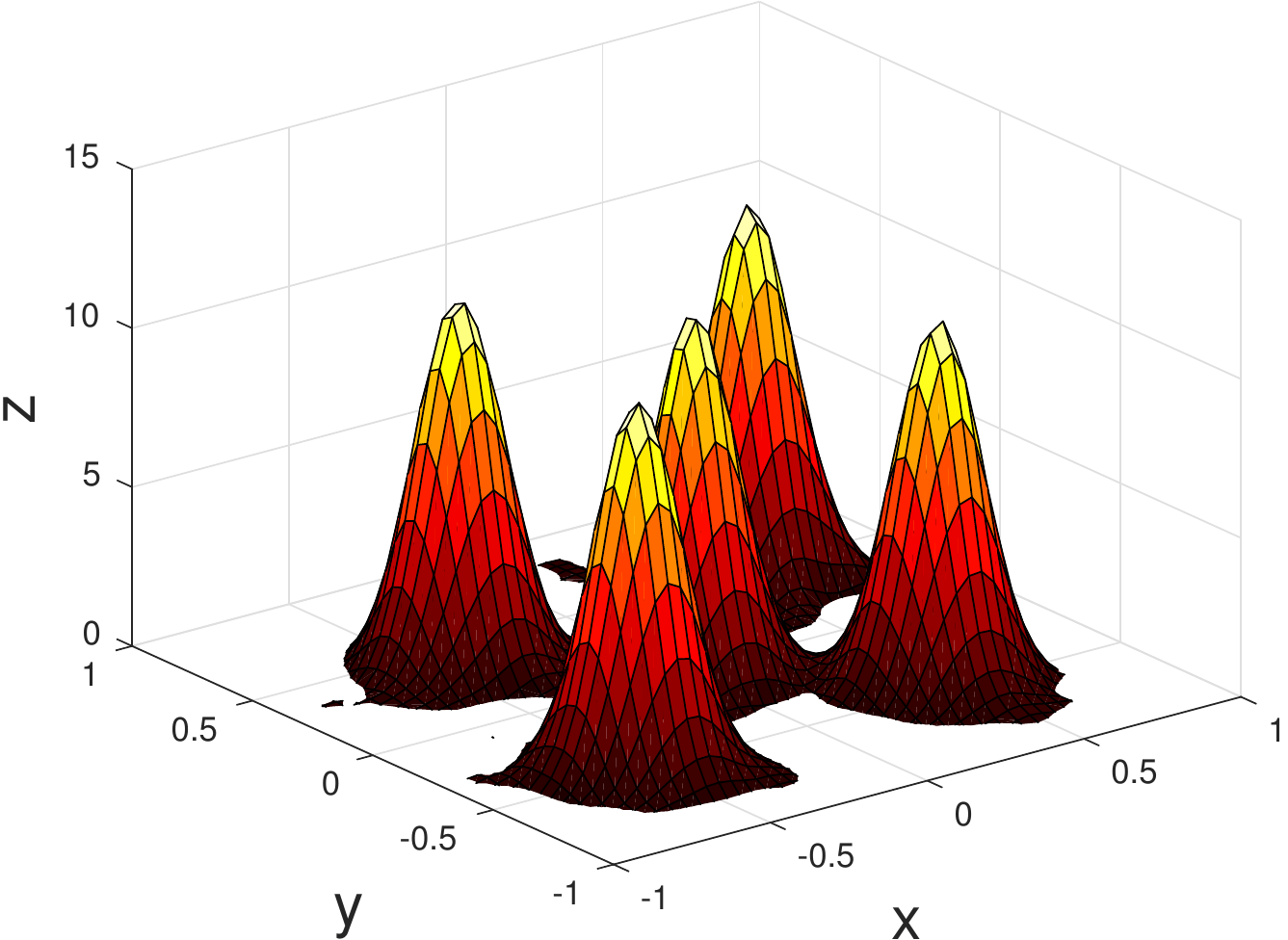}
  \caption{Sample system estimate at $t=1.2$s from a single sensor}
  \label{fig.estimate}
\end{figure}

Consider a physical phenomenon, such as weather, propagating on a planar surface according to the partial differential equation 
\begin{align*}
  \dot u = \left( \alpha \dfrac{\partial^2 u}{\partial x} + \beta \dfrac{\partial ^2 u}{\partial y}\right) .
\end{align*}
As studied in~\cite{khan2008}, and \cite{Kutz2013}, by discretizing the surface into an equally-spaced $n\times n$ grid in the $x$ and $y$ directions and the temporal dimension, thus obtaining a grid of values $u_{ij}^k$, one can construct the linear system representation 
\begin{align*}
\dot{\vec{u}} = \dfrac{2}{\Delta x^2}\left[ I_n\otimes(\alpha D) + (\beta D)\otimes I_n\right] \vec{u}, ~\vec{u} = \Vect(u_{ij}^k)
\end{align*}
where $D$ is a tridiagonal matrix with $1$'s on the sub and superdiagonal, and $-2$'s on the diagonal.
The discrete-time state matrix $A$ can be found using a matrix exponential, or a finite-order Taylor expansion.
Periodic boundary conditions may be added by setting $D_{1n}=D_{n1}=1$.
Suppose that this phenomenon is being measured according to this discretization.
We use a sensor model such that each node collects data from an $n_l\times n_l$ square portion of the grid $u_{ij}^k$ so that each individual grid point is measured by at least one sensor.
The sensor network by the graph in Figure~\ref{fig.graph} is used, with $n_l=4$ for a $50\times50$ grid, and the sensor noise is given unit covariance. 
Five 2D Gaussian functions are used to set the initial conditions for the simulation, which was run over a time interval of 10 seconds.
The true state of the system at 1.2s is shown in Figure~\ref{fig.truestate}, and an aggregate measurement of the entire system is shown in Figure~\ref{fig.ytilde}.
An estimate from a single sensor at this timestep is shown in Figure~\ref{fig.estimate}, and one can see that there is a substantial reduction of noise.

\section{Conclusion}\label{sec.conc}

In this paper, we extended the results of \cite{Sutton1992} to produce a gradient descent estimate of the error covariance in a standard state estimation algorithm
for linear systems.
This gradient descent method was improved with Nesterov acceleration and an adaptive learning rate algorithm.

Finally, the algorithm was extended to improve the Kalman-consensus filter on a sensor network observing a distributed high-dimensional process.
The diagonal estimate of $P$ removes the necessity of computing a large matrix inverse at each timestep at each sensor node, and allows for a computationally cheaper and therefore faster state estimation of the system over the entire network.
The modified Kalman-consensus filter was implemented on a simple 2D diffusion model with 2500 grid points corresponding to the states of the system, and the algorithm was shown to remove a substantial amount of noise from the measurements, and to propagate the state estimate quickly across the sensor network.



\appendix
\section{Gradient Descent for Kalman Filter Covariance Update}\label{app.lin}

Sutton's Gradient Descent for updating the error covariance matrix approximation $\hat P = \mathrm{diag}\left( \hat P_{k}^{(ii)}\right) = \mathrm{diag}\left( \exp \beta^{(i)}_{k+1}\right)$ uses gradient descent on $\beta_{k+1}^{(i)}$: 
\begin{align*}
  \beta_{k+1}^{(i)} = \beta_k^{(i)} - \dfrac{1}{2} \mu \dfrac{\partial \left(\delta_k^T\delta_k\right)}{\partial \beta_k^{(i)}},
\end{align*}
where $\delta_k = \tilde y_k - C \hat x_k^-$ is the error: 
To derive the exact form of the gradient descent on $\beta_k^{(i)}$, consider the following terms: 
\begin{align*}
  \dfrac{\partial \left(\delta^T_k \delta_k\right)}{\partial \hat x_k^{(i)}} &= 2\delta_k^T \xpartial{\delta}{\hat x^{(i)}}{k}{k}= -2\delta_k^T {C^{[i]}}.
\end{align*}
Next, we have 
\begin{align*}
  \xpartial{\delta}{\beta^{(i)}}{k}{k} = \xpartial{~}{\beta^{(i)}}{}{k}\left[ \tilde y_k - C\hat x_k^- \right]
 = -\sum_{j=1}^n C^{[j]} \xpartial{\hat x^{-(j)}}{\beta^{(i)}}{k}{k}.
\end{align*}
We assume that for $i\neq j$,
\begin{align*}
  C^{[j]} \xpartial{\hat x^{-(j)}}{\beta^{(i)}}{k}{k} \cong \vec{0} \implies  \xpartial{\delta}{\beta^{(i)}}{k}{k}  \cong- C^{[i]}   \xpartial{\hat x^{-(i)}}{\beta^{(i)}}{k}{k}.
\end{align*}
Next, define the following two quantities: 
\begin{align*}
  D_k &:= R_k + C\hat P_k C^T = R_k + \sum_{i=1}^n \exp\left({\beta_{k+1}^{(i)}}\right) C^{[i]}C^{[i]T}\\
  k_k^{(i)}&:= \exp\left({\beta_{k+1}^{(i)}}\right) D^{-1}_k C^{[i]}.
\end{align*}
Next, we compute the following derivatives: 
\begin{align*}
  \xpartial{D}{\beta^{(j)}}{k}{} &= \sum_{i=1}^n \dfrac{\partial e^{\beta_{k+1}^{(i)}}}{\partial\beta^{(j)}}C^{[i]}C^{[i]T} = { e_{k+1}^{\beta^{(j)}}}C^{[j]}C^{[j]T}\\
 \dfrac{\partial D^{-1}_k}{\partial \beta^{(j)}} &=  -D_k^{-1}  \dfrac{\partial D_k}{\partial \beta^{(j)}} D_k^{-1}= -k_k^{(i)}C^{[j]T}D_k^{-1}\\
\dfrac{\partial k_{k+1}^{(i)}}{\partial \beta^{(i)}} &= \left[ \dfrac{\partial e^{\beta_{k+1}^{(i)}}}{\partial \beta ^{(i)}}  D^{-1}_{k} + e^{ \beta ^{(i)}_{k+1}}\dfrac{ \partial D_k^{-1}  }{\partial \beta^{(i)}}\right]C^{[i]}\\
 &=\beta_{k+1}^{(i)} D_k^{-1} C^{[i]} -  e^{\beta ^{(i)}_{k+1}}k_{k+1}^{(i)}C^{[i]T}D_k^{-1}C^{[i]}\\
&= k_{k+1}^{(i)} - k_{k+1}^{(i)}\left( e^{\beta ^{(i)}_{k+1}} \left(D_k^{-1}\right)^T C^{[i]}\right)^TC^{[i]}\\
&= k_{k+1}^{(i)}- k_{k+1}^{(i)}k_{k+1}^{(i)T}C^{[i]}.
\end{align*}
Now, we can write the gradient descent on $\beta^{(i)}$: 
\begin{align*}
  \beta_{k+1}^{(i)} &= \beta_k^{(i)} - \dfrac{1}{2}\mu \dfrac{\partial (\delta_k^T\delta)}{\partial \beta^{(i)}}\\
  &= \beta_k^{(i)} - \dfrac{1}{2}\mu \sum_j \dfrac{\partial (\delta_k^T\delta)}{\partial \hat x _k^{-(j)}}\dfrac{\partial\hat x_k^{-(j)}}{\partial \beta ^{(i)}}\\
  &\cong \beta_k^{(i)} - \dfrac{1}{2}\mu \dfrac{\partial (\delta_k^T\delta)}{\partial \hat x _k^{-(i)}}\dfrac{\partial\hat x_k^{-(i)}}{\partial \beta ^{(i)}}\\
  &= \beta_k^{(i)} + \mu \delta_k^T C^{[i]} \dfrac{\partial\hat x_k^{-(i)}}{\partial \beta^{(i)}}.
\end{align*}
All that remains is to find an approximation to $\frac{\partial\hat x_k^{-(i)}}{\partial \beta^{(i)}}$.
We denote this approximation with $h_{k+1}^{(i)}$:
\begin{align*}
  h_{k+1}^{(i)} &\cong \dfrac{\partial\hat x_{k+1}^{-(i)}}{\partial \beta^{(i)}} = \dfrac{\partial}{\partial\beta^{(i)}} \left[ \hat x_k^{(i)} + k_k^{(i)T}\delta_k\right]\\
&= \dfrac{\partial\hat x_k^{(i)}}{\partial \beta^{(i)}} + \dfrac{\partial k_k^{(i)T}}{\partial \beta^{(i)}} \delta_k + k_k^{(i)T}\dfrac{\partial \delta_k}{\partial \beta^{(i)}}\\
&= \dfrac{\partial\hat x_k^{(i)}}{\partial \beta^{(i)}} + \left( k_k^{(i)} - k_k^{(i)}k_k^{(i)T} C^{[i]}\right)^T \delta_k - k_k^{(i)T}\dfrac{\partial \hat x_k^{(i)}}{\partial \beta^{(i)}} C^{[i]}\\[-3.5ex]
&= h_k^{(i)}\left(1 - k_k^{(i)T}C^{[i]}\right) + \left( k_k^{(i)} - k_k^{(i)}k_k^{(i)T} C^{[i]}\right)^T \delta_k.
\end{align*}
This yields the final gradient descent update: 
\begin{align*}
  \beta_{k+1}^{(i)} = \beta_k^{(i)} + \mu \delta^T_k C^{[i]}h_k^{(i)}.
\end{align*}
In conclusion, the Kalman filter with gradient descent error covariance update is given by the equations
\begin{align*}
 \hat x_{k+1} &= \hat x_k + K_k \left[ \tilde y_k - C \hat{x}_k\right] \\
  K_k &= \hat P_k C^T \left[ C \hat P_k C^T + R_k\right]^{-1}\\
\hat P &=  \mathrm{diag}\left( e ^{\beta^{(i)}_{k+1}}\right)\\
  \beta_{k+1}^{(i)} &= \beta_k^{(i)} + \mu \delta^T_k C^{[i]}h_k^{(i)}\\
h_{k+1}^{(i)} &= h_k^{(i)}\left(1 - k_k^{(i)T}C^{[i]}\right)^+ + \left( k_k^{(i)} - k_k^{(i)}k_k^{(i)T} C^{[i]}\right)^T \delta_k.
\end{align*}


\end{document}